\documentclass[twoside]{amsart}

\usepackage{latexsym}
\usepackage{amsthm,amsfonts,amssymb}
\usepackage[leqno]{amsmath}
\usepackage[all]{xy} \SelectTips{eu}{}
\usepackage{hyperref}

\newcommand{\numberseries}{\mdseries}   

\newlength{\thmtopspace}                
\newlength{\thmbotspace}                
\newlength{\thmheadspace}               
\newlength{\thmindent}                  

\newtheoremstyle{bfupright head,slanted body}
                {\thmtopspace}{\thmbotspace}
                {\slshape}{\thmindent}{\bfseries}{.}{\thmheadspace}
                {{\numberseries \thmnumber{(#2) }}\thmnote{#3}}

\newtheoremstyle{bfupright head,upright body}
                {\thmtopspace}{\thmbotspace}
                {\upshape}{\thmindent}{\bfseries}{.}{\thmheadspace}
                {{\numberseries \thmnumber{(#2) }}\thmnote{#3}}

\newtheoremstyle{bfit head,upright body}
                {\thmbotspace}{\thmbotspace}
                {\upshape}{\thmindent}{\upshape}{.}{\thmheadspace}
                {{\numberseries\thmnumber{(#2) }}
                {\bfseries\itshape\thmnote{\negthickspace#3}}}

\newtheoremstyle{fixed bf head,slanted body}
                {\thmtopspace}{\thmbotspace}{\slshape}
                {\thmindent}{\bfseries}{.}{\thmheadspace}
                {{\numberseries \thmnumber{(#2) }}\thmname{#1}\thmnote{ (#3)}}

\newtheoremstyle{fixed bf head,upright body}
                {\thmtopspace}{\thmbotspace}{\upshape}
                {\thmindent}{\bfseries}{.}{\thmheadspace}
                {{\numberseries \thmnumber{(#2)
                    }}\thmname{#1}\thmnote{ (#3)}}

\newtheoremstyle{numbered paragraph}
                {\thmtopspace}{\thmbotspace}{\upshape}
                {\thmindent}{\upshape}{}{0pt}
                {{\numberseries \thmnumber{(#2) }}}

\newtheoremstyle{independent paragraph}
                {\thmtopspace}{\thmbotspace}
                {\upshape}{\parindent}{\upshape}{}{0pt}
                {\thmnote{#3 }}

\newtheoremstyle{subparagraph}
                {\thmbotspace}{\thmbotspace}
                {\upshape}{\parindent}{\upshape}{}{0pt}
                {\thmnote{#3 }}

\newtheoremstyle{notes}
                {\thmtopspace}{\thmbotspace}
                {\ttfamily}{\thmindent}{\ttfamily\small }{}{0pt}
                {\thmnote{#3 }}

\theoremstyle{bfupright head,slanted body}
\newtheorem{res}{}             \newtheorem*{res*}{}

\theoremstyle{bfupright head,upright body}
\newtheorem{bfhpg}[res]{}               \newtheorem*{bfhpg*}{}

\theoremstyle{bfit head,upright body}

\theoremstyle{fixed bf head,slanted body}
\newtheorem{thm}[res]{Theorem}          \newtheorem*{thm*}{Theorem}
      \newtheorem*{prp*}{Proposition}
        \newtheorem*{cor*}{Corollary}
\newtheorem{lem}[res]{Lemma}            \newtheorem*{lem*}{Lemma}

\theoremstyle{fixed bf head,upright body}
       \newtheorem*{dfn*}{Definition}
      \newtheorem*{obs*}{Observation}
           \newtheorem*{rmk*}{Remark}
\newtheorem{exa}[res]{Example}          \newtheorem*{exa*}{Example}

\theoremstyle{numbered paragraph}

\newlength{\thmlistleft}        
\newlength{\thmlistright}       
\newlength{\thmlistpartopsep}   
\newlength{\thmlisttopsep}      
\newlength{\thmlistparsep}      
\newlength{\thmlistitemsep}     

\newcounter{eqc} 
  {\end{list}}%

\newcounter{prt}
  {\end{list}}%

\newcounter{rqm}
  {\end{list}}%

  {\end{list}}%

  {\end{list}}%

\newenvironment{prf}[1][Proof]{%
\begin{proof}[\bf #1]
\setcounter{equation}{0}
\renewcommand{\theequation}{\arabic{equation}}}
{\end{proof}}

\renewcommand{\eqref}[1]{\pgref{eq:#1}}

\newcommand{\pgref}[1]{(\ref{#1})}

\newcommand{\lemref}[2][Lemma~]{#1\pgref{lem:#2}}
\newcommand{\thmref}[2][Theorem~]{#1\pgref{thm:#2}}
\newcommand{\corcite}[2][?]{\cite[cor.~#1]{#2}}

\newcommand{\prpcite}[2][?]{\cite[prop.~#1]{#2}}
\newcommand{\thmcite}[2][?]{\cite[thm.~#1]{#2}}

\newcommand{\mgen}[1]{\nu{(#1)}}
\newcommand{\s}[1][R]{\mathfrak{S}_0(#1)}
\newcommand{\tp}[3][R]{#2\mspace{-1mu}\otimes_{\mspace{-1mu}#1}\mspace{-1mu}#3}
\newcommand{\Hom}[3][R]{\operatorname{Hom}_{#1}(#2,#3)}
\newcommand{\Ext}[4][R]{\operatorname{Ext}_{#1}^{#2}(#3,#4)}
\newcommand{\no}{\mspace{-1mu}} 
\newcommand{\f}{\varphi}
\newcommand{\is}{\cong}
\newcommand{\x}{\pmb{x}}
\newcommand{\Rhat}{\widehat{R}}
\newcommand{\onto}{\twoheadrightarrow}
\newcommand{\mapdef}[4][\rightarrow]{\mbox{\ensuremath{#2\colon #3 #1 #4}}}
\newcommand{\lgt}[2][R]{\operatorname{length}_{#1}#2}

\renewcommand{\theequation}{\arabic{equation}}
\numberwithin{equation}{res}

\begin{document}

\title[A C.-M.\ algebra has only finitely many semidualizing
  modules]
{A Cohen-Macaulay algebra has only finitely many semidualizing
  modules}

\author{Lars Winther Christensen}

\address{Department of Math.\ and Stat., Texas Tech University,
  Lubbock, TX 79409, U.S.A.}

\email{lars.w.christensen@ttu.edu}

\urladdr{http://www.math.ttu.edu/$\scriptstyle\sim$lchriste}

\thanks{This work was done while L.W.C.\ visited University of
  Nebraska--Lincoln, partly supported by grants from the Danish
  Natural Science Research Council and the Carlsberg Foundation.}

\author{Sean Sather-Wagstaff}

\address{Department of Mathematics,
300 Minard Hall,
North Dakota State University,
Fargo, ND 58105-5075, U.S.A.}

\email{Sean.Sather-Wagstaff@ndsu.edu}

\urladdr{http://math.ndsu.nodak.edu/faculty/ssatherw/}

\thanks{S.S.-W.\ was partially supported by NSF grant NSF~0354281.}


\date{16 November 2007}

\keywords{Selforthogonal, self-orthogonal, semidualizing,
  semi-dualizing, spherical modules}

\subjclass[2000]{13C13, 13H10}


\begin{abstract}
  We prove the result stated in the title, which answers the
  equicharacteristic case of a question of Vasconcelos.
\end{abstract}

\maketitle

\vspace*{\stretch{1}}

In this paper, $R$ is a commutative noetherian local ring.  A finitely
generated $R$-module $C$ is \emph{semidualizing} if the homothety
morphism $R\to\Hom{C}{C}$ is an isomorphism and
$\Ext[R]{\ge1}{C}{C}=0$. Examples include $R$ itself and (if one
exists) a dualizing $R$-module in the sense of Grothendieck.
Semidualizing modules have properties similar to those of dualizing
modules \cite{LWC01a} and arise in several contexts.

Let $\s$ denote the set of isomorphism classes of semidualizing
$R$-modules.  Vasconcelos, calling these modules ``spherical,'' asked
whether $\s$ is finite when $R$ is Cohen-Macaulay and whether it has
even cardinality when it contains more than one
element~\cite[p.~97]{dtmc}. In \cite{SSW} affirmative answers to these
questions are given, e.g., for certain determinantal rings. Here we
prove:

\begin{thm}
  \label{thm:even}
  If $R$ is Cohen-Macaulay and equicharacteristic, then $\s$ is
  finite.
\end{thm}

\begin{bfhpg}[Remark]
This result also yields an answer to the parity  part of Vasconcelos' question
for certain Cohen-Macaulay rings. 
Let $R$ be as in \thmref{even}.
If $R$ is Gorenstein, then $\s$ has exactly one element, namely the isomorphism
class of $R$; see 
\corcite[(8.6)]{LWC01a}.
On the other hand, if $R$ has a dualizing module and is not Gorenstein, then
$\s$ has even cardinality. Indeed, in the 
derived category of $R$, every semidualizing $R$-complex 
is isomorphic to a shift of a semidualizing $R$-module; see \corcite[3.4]{AFrSSW07b}. 
Thus, the set of shift-isomorphism classes of 
semidualizing $R$-complexes is finite,
and so
\prpcite[(3.7)]{LWC01a} shows that  it has even cardinality.
\end{bfhpg}

In preparation for the proof of \thmref{even}, we recall some recent results
that allow us to reduce it to a problem  in representation theory
over artinian $k$-algebras.

\begin{bfhpg}[Fact]
  \label{tp}
  Let $\f\colon R \to S$ be a local ring homomorphism of finite flat
  dimension, i.e., such that $S$ has finite flat dimension as an
  $R$-module via $\f$. If $C$ is a semidualizing $R$-module, then
  $\tp{S}{C}$ is a semidualizing $S$-module
  by~\prpcite[(5.7)]{LWC01a}.  If $C$ and $C'$ are semidualizing
  $R$-modules such that $\tp{S}{C}$ and $\tp{S}{C'}$ are isomorphic as
  $S$-modules, then $C$ and $C'$ are isomorphic as $R$-modules; see
  \thmcite[4.5 and 4.9]{AFrSSW07a}.  Thus, the functor $\tp{S}{-}$
  induces an injective map $\s[\f]\colon\s \hookrightarrow \s[S]$.
\end{bfhpg}

\begin{lem}
  \label{lem:finite}
  Assume there are local ring homomorphisms of finite flat dimension
  \begin{equation*}
    \xymatrix{R \ar[r]^-{\f} & R' & Q \ar@{>>}[l]_-{\mspace{20mu} \rho}
      \ar[r]^-{\tau} & Q'}
  \end{equation*}
  such that $\rho$ is surjective with kernel generated by a
  $Q$-regular sequence.  Then there are inequalities of cardinalities $|\s|
  \leq |\s[\widehat{Q}]| \leq |\s[\widehat{Q'}]|$.
\end{lem}

\begin{prf}
  The completion morphism $\epsilon\colon R\to\Rhat$ 
  conspires with the completions of the given maps to yield the following
  \begin{equation*}
    \xymatrix{
      \s \ar@{^{(}->}[rr]^-{\s[\hat{\f}\epsilon]} 
      && \mathfrak{S}_0\bigl(\widehat{R'}\bigr)
      && \mathfrak{S}_0\bigl(\widehat{Q}\bigr) 
      \ar@{_{(}->>}[ll]_-{\s[\hat{\rho}]} 
      \ar@{^{(}->}[rr]^-{\s[\hat{\tau}]} 
      && \s[\widehat{Q'}]. }
  \end{equation*}
  The injectivity of the induced maps is justified in Fact \pgref{tp};
  the surjectivity of $\s[\hat{\rho}]$ follows from
  \prpcite[4.2]{AFrSSW07a} because $\widehat{Q}$ is complete and
  $\hat{\rho}$ is surjective with kernel generated by a $\widehat{Q}$-regular
  sequence.  The desired inequalities now follow.
\end{prf}

\begin{bfhpg}[Proof of (1)]
  Let $\x$ be a system of parameters for $R$ and set $R'=R/(\x)$ with
  $\mapdef{\varphi}{R}{R'}$ the natural surjection.  Using
  \lemref{finite} with $Q'=Q=R' \is \widehat{R'}$, we may replace $R$
  with $R'$ in order to assume that $R$ is artinian.
  
  There is a flat homomorphism of artinian local rings $R \to R''$,
  such that $R''$ has algebraically closed residue field; see
  \prpcite[0.(10.3.1)]{egaIII1}. By Fact \pgref{tp} we may replace $R$
  by $R''$ to assume that its residue field $k$ is algebraically
  closed.  As $R$ is equicharacteristic and artinian, Cohen's
  structure theorem implies $R$ is a $k$-algebra.
  
  For an $R$-module $M$, let $\mgen{M}$ denote the minimal number of
  generators of~$M$. Let $C$ be a semidualizing $R$-module and let $E$
  be the injective hull of $k$.  Hom-evaluation \prpcite[5.3]{careil}
  and the homothety map yield a sequence of isomorphisms
  $$
  \tp{C}{\Hom{C}{E}}
  \is \Hom{\Hom{C}{C}}{E}
  \is\Hom{R}{E}
  \is E.$$
  Hence, there is an inequality $\mgen{C} \le
  \mgen{E}$. This gives the second inequality below; the first is
  from a surjection $R^{\mgen{C}} \onto C$.%
  \begin{equation*}
    \lgt{C} \le \mgen{C}\cdot\lgt[\no]{R} \le \mgen{E}\cdot\lgt[\no]{R}
  \end{equation*}
  By \cite[proof of first prop.\ in sec.~3]{DHp95} there are only
  finitely many isomorphism classes of $R$-modules $M$ with $\Ext{\ge
    1}{M}{M}=0$ and \mbox{$\lgt{M} \le \mgen{E}\cdot\lgt[\no]{R}$}.
  From the displayed inequalities it follows that $\s$ is a finite
  set.  \qed
\end{bfhpg}

Finally we illustrate how \thmref[]{even} and \lemref[]{finite} apply
to answer Vasconcelos' finiteness question for certain rings of mixed
characteristic.

\begin{exa}
  Let $Q$ be a complete Cohen-Macaulay local ring with residue field
  of characteristic $p>0$.  If $p$ is $Q$-regular, then $\s$ is finite
  for every local ring $R$ such that $\widehat{R}\cong Q$ or
  $\widehat{R}\cong Q/(p^n)$ for some $n\geq 2$: \thmref{even} implies
  that $\s[Q/(p)]$ is finite, and \lemref{finite} applies to the
  diagram $R\to\Rhat\leftarrow Q\to Q/(p)$.
\end{exa}


\section*{Acknowledgments}
\label{sec:ack}

We thank Ragnar-Olaf Buchweitz, Hubert Flenner, Graham Leuschke, and
Diana White for useful conversations related to this work.  We also
thank the referee for thoughtful suggestions.


\bibliographystyle{amsplain}
\def\cprime{$'$} \providecommand{\bysame}{\leavevmode\hbox
  to3em{\hrulefill}\thinspace}
\providecommand{\MR}{\relax\ifhmode\unskip\space\fi MR }
\providecommand{\MRhref}[2]{%
  \href{http://www.ams.org/mathscinet-getitem?mr=#1}{#2} }
\providecommand{\href}[2]{#2}

\end{document}